\documentclass{article}

\usepackage{amssymb}
\usepackage{latexsym}
\usepackage{amsthm}
\usepackage{enumerate}
\usepackage{epsfig}
\usepackage{color}
\usepackage{float}
\usepackage{subfigure}
\usepackage{amsmath}
\usepackage{ifthen}
\usepackage{lscape}
\usepackage{setspace}
\usepackage{amscd}
\usepackage{xspace}
\usepackage[nobysame,alphabetic]{amsrefs}
\usepackage{srcltx}

\def\a{\alpha}
\def\b{\beta}
\def\C{\mathcal{C}}
\def\g{\gamma}
\def\G{\Gamma}

\def\F{\mathcal{F}}

\def\L{\Lambda}
\def\l{\lambda}

\def\O{\mathcal{O}}
\def\R{\mathbb{R}}

\def\S{\Sigma}
\def\T{\mathcal{T}}
\def\Z{\mathbb{Z}}
\def\d{\partial}
\def\cross{\times}

\def\e{\epsilon}

\def\t{\tau}

\def\le{\leqslant}
\def\ge{\geqslant}

\def\MF{\mathcal{MF}}

\def\PMF{\mathcal{PMF}}
\def\MF{\mathcal{MF}}

\def\Grel{\widehat G}
\def\m{{\bf m}}
\def\Re{\text{Re}}
\def\Im{\text{Im}}

\def\etabar{\overline{\eta}}

\def\fmin{\mathcal{F}_{min}}

\def\Teich{Teichm\"uller }

\def\mcg{mapping class group\xspace}
\def\pA{pseudo-Anosov\xspace}
\def\cone{\text{Cone}}

\def\abem{Athreya, Bufetov, Eskin and Mirzakhani\xspace}

\newcommand{\norm}[1]{|#1|}
\newcommand{\nhat}[1]{|\widehat{#1}|}
\newcommand{\dhat}[1]{\widehat d (#1)}
\newcommand{\dc}[1]{d_\C(#1)}
\newcommand{\dt}[1]{d_\T(#1)}

\newcommand{\horo}[2]{\widehat{\O}_{#2}(#1)}
\newcommand{\sect}[2]{Sect_{#1}(#2)}
\newcommand{\gp}[2]{(#1|#2)}
\newcommand{\conj}[1]{\nhat{#1}_c}

\newtheorem{theorem}{Theorem}[section]
\newtheorem{lemma}[theorem]{Lemma}

\newtheorem{proposition}[theorem]{Proposition}

\theoremstyle{definition}

\newtheorem{claim}[theorem]{Claim}

\restylefloat{figure}

\begin{document}


\title{Asymptotics for pseudo-Anosov elements in \Teich lattices}
\author{Joseph Maher\footnote{joseph.maher@csi.cuny.edu}}
\date{\today}

\maketitle

\begin{abstract}
A \emph{\Teich lattice} is the orbit of a point in \Teich space under
the action of the mapping class group. We show that the proportion of
lattice points in a ball of radius $r$ which are not \pA tends to zero
as $r$ tends to infinity. In fact, we show that if $R$ is a subset of the
\mcg, whose elements have an upper bound on their translation length on
the complex of curves, then proportion of lattice points in the ball
of radius $r$ which lie in $R$ tends to zero as $r$ tends to infinity.

Subject code: 37E30, 20F65, 57M50.

\end{abstract}

\tableofcontents

\section{Introduction}

A \emph{\Teich lattice} $\G y$ is the orbit of a point $y$ in \Teich
space under the action of the mapping class group $\G$.  Athreya,
Bufetov, Eskin and Mirzakhani \cite{abem} showed that the asymptotic
growth rate of the number of lattice points in a ball of radius $r$ is
\[ \norm{\G y \cap B_r(x)} \sim \L(x) \L(y) h e^{hr}. \]
Here $B_r(x)$ denotes the ball of radius $r$ centered at $x$ in the
\Teich metric, $h = 6g-6$ is the topological entropy of the \Teich
geodesic flow, $\L$ is the Hubbard-Masur function, $\norm{X}$ denotes
the number of elements in a finite set $X$, and $f(r) \sim g(r)$ means
$f(r)/g(r)$ tends to one as $r$ tends to infinity. We use their work,
together with some results from \cite{maher1}, to show that the number
of lattice points corresponding to \pA elements in the ball of radius
$r$ is asymptotically the same as the total number of lattice points
in the ball of radius $r$. More generally, let $R \subset \G$ be a set
of elements for which there is an upper bound on their translation
length on the complex of curves, for example, the set of non-\pA
elements.  We shall write $Ry$ for the orbit of the point $y$ under
the subset $R$. We show that the proportion of lattice points $\G y$
in $B_r(x)$ which lie in $R y$ tends to zero as $r$ tends to infinity.
In fact, we show a version of this result for bisectors.  Let $Q$ be
the space of unit area quadratic differentials on the surface $\S$,
and given $x \in \T$, let $S(x)$ be the subset of $Q$ consisting of
unit area quadratic differentials on $x$. The space $Q$ has a
canonical measure, known as the Masur-Veech measure, which we shall
denote $\mu$, and we will write $s_x$ for the conditional measure on
$S(x)$ induced by $\mu$.  We may think of $S(x)$ as the (co-)tangent
space at $x$.  Given $x, y \in \T$, let $q_x(y)$ be the unit area
quadratic differential on $x$ corresponding to the geodesic ray
through $y$.  Given a lattice point $\g y$, we will write $q(x, \g y)$
for the pair $(q_x(\g y), q_y(\g^{-1} x)) \in S(x) \cross S(y)$.
Given subsets $U \subset S(x)$ and $V \subset S(y)$, we may consider
those lattice points $\g y$ which lie in the bisector determined by
$U$ and $V$, i.e. those $\g y$ for which $q(x, \g y) \in U \cross V$.
If $X$ is a finite subset of $\G$, we will write $\norm{X,
  \text{\emph{condition}}}$ to denote the number of elements $\g \in
X$ which also satisfy \emph{condition}. We say a surface of finite
type is \emph{sporadic} if it is a torus with at most one puncture, or
a sphere with at most four punctures.

\begin{theorem} \label{theorem:main}
Let $\G$ be the \mcg of a non-sporadic surface.
Let $R \subset \G$ be a set of elements of the for which there is an upper
bound on their translation distance on the complex of curves. Let $x,
y \in \T$, and let $U \subset S(x)$ and $V \subset S(y)$ be Borel sets
whose boundaries have measure zero. Then
\begin{equation} \label{eq:main} 
\frac{ \norm{R y \cap B_r(x), \ q(x, \g y) \in U \cross V} }{
  \norm{ \G y \cap B_r(x), \ q(x, \g y) \in U \cross V } } \to 0,
\text{ as } r \to \infty. 
\end{equation}
\end{theorem}

This shows that \pA elements are ``generic'' in the \mcg, at least for
one particular definition of generic, see Farb \cite{farb2} for a
discussion of similar questions.  In the case in which $R$ consists of
the non-\pA elements of the \mcg, this result should also follow from
the methods of Eskin and Mirzakhani \cite{em}, which they use to show
that the number of conjugacy classes of \pA elements of translation
length at most $r$ on \Teich space is asymptotic to $e^{hr}/hr$. In
the sporadic cases, the mapping class group is either finite, or
already well understood, as the mapping class group is $SL(2,\Z)$, up 
to finite index.

In the remainder of this section we give a brief outline of the
argument. In Section \ref{section:geodesic flow} we describe the
results we need from \abem \cite{abem} and set up some notation. In
Section \ref{section:limit sets} we review some useful properties of
the visual boundary of \Teich space, and then in Section
\ref{section:asymptotics} we prove the main result.

\subsection{Outline} 

Let $R \subset \G$ be a set of elements for which there is an upper
bound on their translation length on the complex of curves, for
example, the set of non-\pA elements in the \mcg. We wish to consider
the distribution of elements of $Ry$ inside \Teich space $\T$. In some
parts of $\T$ elements of $Ry$ are close together, and in other parts
elements of $Ry$ are widely separated. We quantify this by by defining
$R_k$ to be the $k$-dense elements of $Ry$, namely those elements of
$Ry$ which are distance at most $k$ in the \Teich metric from some
other element of $R$. If two lattice points $\g y$ and $\g' y$ are
a bounded \Teich distance apart, then $\g$ and $\g'$ are a bounded
distance apart in the word metric on $\G$.  In \cite{maher1} we showed
that the limit set of the $k$-dense elements in the word metric has
measure zero in the Gromov boundary of the relative space, and we use
this to show that the $k$-dense elements in \Teich space have a limit
set in the visual boundary $S(x)$ which has $s_x$-measure zero. A
straightforward application of the results of \cite{abem} then shows
that the proportion of lattice points $\G y \cap B_r(x)$ which lie in
$R_k y$ tends to zero as $r$ tends to infinity.

It remains to consider $R y \setminus R_k y$, which we shall denote
$R'_k y$. We say a subset of $\T$ is $k$-separated, if any two
elements of the set are \Teich distance at least $k$ apart, so $R'_k
y$ is a $k$-separated subset of $\T$.  Naively, one might hope that
the proportion of $k$-separated elements of $\G y$ in $B_r(x)$ is at
most $1 / \norm{\G y \cap B_k(y)}$, as each lattice point $\g y \in
R'_k y$ is contained in a ball of radius $k$ in \Teich space
containing $\norm{\G y \cap B_k(y)}$ other lattice points, none of
which lie in $R'_k y$. Such a bound would imply the required result,
as this would give a collection of upper bounds for
\[ \lim_{r \to \infty} \frac{ \norm{ R y \cap B_r(x)} }{\norm{\G y
    \cap B_r(x)}} \]
which depend on $k$, and furthermore these upper bounds would decay
exponentially in $k$, so this implies that the limit above is zero.
However, this argument only works for those lattice points in the
interior of $B_r(x)$ for which $B_k(\g y) \subset B_r(x)$.  If a
lattice point $\g y$ is within distance $k$ of $\d B_r(x)$, then many
of the lattice points in $B_k(\g y)$ may lie outside $B_r(x)$, and a
definite proportion of lattice points are close to the boundary, as
the volume of $B_r(x)$ grows exponentially. We use the mixing property
of the geodesic flow to show that that $\d B_r(x)$ becomes
equidistributed on compact sets of the quotient space $\T / \G$, and
this in turn shows that the intersections of $\d B_r(x)$ with $B_k(\g
y)$ are evenly distributed.  This implies that we can find an upper
bound for the average number of lattice points of $\G y$ near some $\g
y$ close to the boundary, which do in fact lie inside $B_r(x)$.  In
fact, we prove a result that works for bisectors, so we also need to
show that the proportion of lattice points near the geodesics rays
through $\d U$ tends to zero as $r$ tends to infinity.  These
arguments using mixing originate in work of Margulis \cite{margulis},
and our treatment of conditional mixing is essentially due to Eskin
and McMullen \cite{emc}, see also Gorodnik and Oh \cite{go}, for the
higher rank case.

\subsection{Acknowledgements} 

I would like to thank Alex Eskin, Howard Masur and Kasra Rafi for
useful advice. This work was partially supported by NSF grant
DMS-0706764.

\section{The \Teich geodesic flow} \label{section:geodesic flow}

In this section we review the work of Athreya, Bufetov, Eskin and
Mirzakhani \cite{abem} that we will use, fix notation, and use the
mixing property of the geodesic flow to show a conditional mixing
result, which is an analogue in \Teich space of a result of Eskin and
McMullen \cite{emc} in the case of Lie groups.

Let $\S_{g,b}$ be an orientable surface of finite type, of genus $g$
and with $b$ punctures, which is not a torus with one or fewer
punctures, or a sphere with four or fewer punctures. We will just
write $\S$ for the surface if we do not need to explicitly refer to
the genus or number of punctures. Let $\G$ be the mapping class group
of $\S$, and we will consider $\G$ to be a metric space with the word
metric coming from some fixed choice of generating set.  We will write
$\T$ for the \Teich space of conformal structures on $\S$, with the
\Teich metric, and $\T$ is homeomorphic to $\R^{6g-6+2b}$. We will
write $B_r(x)$ for the ball of radius $r$, centered at $x$ in $\T$. A
choice of basepoint $y$ for $\T$ determines a map from $\G$ to $\T$,
defined by $\g \mapsto \g y$.  We shall write $\G y$ for the image of
$\G$ under this map, which we shall call a \Teich lattice. The map $\g
\mapsto \g y$ is coarsely distance decreasing, but is not a
quasi-isometry.

Let $\MF$ be the space of measured foliations on the surface $\S$,
which is homeomorphic to $\R^{6g-6+2b}$, and let $\nu$ be the Thurston
measure on $\MF$, which is preserved by the action of $\G$.  Let $Q$
be the space of unit-area quadratic differentials on $\S$, let $\pi:Q
\to \T$ be the projection from a quadratic differential to the
underlying Riemann surface, and let $S(x)$ be the pre-image of $x$ in
$Q$ under the projection.  Hubbard and Masur \cite{hm} showed that for
$x \in \T$, the map which sends a quadratic differential $q$ on $x$ to
its vertical foliation $\Re(q^{1/2})$ is a homeomorphism, and we shall
write $\eta^+ : Q \to \MF$ for the restriction of this map to unit
area quadratic differentials, and $\overline{\eta}^+: Q \to \PMF$ for
the induced map from unit area quadratic differentials to projective
equivalence classes of their vertical foliations.  Similarly, the map
that sends a quadratic differential to its horizontal foliation
$\Im(q^{1/2})$ is a homeomorphism, and we shall write $\eta^- : Q \to
\MF$ for the restriction of this map to unit area quadratic
differentials, and $\etabar^-$ for the corresponding map to $\PMF$. In
particular, the restriction $\etabar^+ : S(x) \to \PMF$ is a
homeomorphism, as is the restriction of $\etabar^-$.

Masur \cite{masur} and Veech \cite{veech2} showed that the space $Q$
carries a $\G$-invariant smooth measure $\mu$, preserved by the \Teich
geodesic flow, such that $\mu(Q/\G)$ is finite, and this measure is
unique up to rescaling, so we shall assume that $\mu(Q / \G) = 1$. A
quadratic differential $q$ is uniquely determined by its real and
imaginary measured foliations, $\eta^+(q)$ and $\eta^-(q)$, so the map
$\eta^+ \cross \eta^-$ gives an embedding of $Q$ in $\MF \cross \MF$.
The Masur-Veech measure $\mu$ is then defined by $\mu(E) = (\nu \cross
\nu) (\cone(E))$, where $\cone(E)$ is the cone over $E$ based at the
origin, i.e. $\{ tq \mid q \in E, 0 < t \le 1 \}$.  We shall write
$\m$ for the induced measure on \Teich space, i.e.  $\m = \pi_* \mu$.

Given a point $x \in \T$, the visual boundary at $x$ is the space of
geodesic rays based at $x$, which may be identified with $S(x)$, the
space of unit area quadratic differentials on $x$.  We shall write
$\overline {\T}_x$ for $\T \cup S(x)$, the compactification of \Teich
space using the visual boundary at $x$.  We will write $q_x(y)$ for
the unit area quadratic differential on $x$ which corresponds to the
\Teich geodesic ray starting at $x$ which passes through $y$.  Given a
subset $U$ of $S(x)$ we shall write $\sect{x}{U}$ for the union of
geodesic rays based at $x$ corresponding to quadratic differentials in
$U$, so $y \in \sect{x}{U}$ if and only if $q_x(y) \in U$.

Let $g_t$ be the \Teich geodesic flow on $Q$. The flow $g_t$ commutes
with the action of $\G$, preserves the measure $\mu$, and preserves
the following foliations.  The \emph{strong stable foliation}
$\F^{ss}$ has leaves of the form $\{ q \in Q \mid \eta^+(q) =
\text{const} \}$, and if $p$ is a point in $Q$ we will write
$\a^{ss}(p)$ for the leaf through $p$, i.e.
\[ \a^{ss}(p) = \{ q \in Q \mid \eta^+(q) = \eta^+(p) \}. \]
The \emph{strong unstable foliation} has leaves of the form $\{ q \in Q
\mid \eta^-(q) = \text{const} \}$, and we shall write $\a^{ss}(p)$ for
the leaf through $p$, i.e.
\[ \a^{uu}(p) = \{ q \in Q \mid \eta^-(q) = \eta^-(p) \}. \]
These foliations are $\G$-invariant, so they descend to foliations on
$Q / \G$. The \emph{unstable foliation} $\F^u$ has leaves of the form
\[ \a^u(q) = \bigcup_{t \in \R} g_t \a^{uu}(q),\] 
while the \emph{stable foliation} $\F^s$ has leaves of the form 
\[ \a^s(q) = \bigcup_{t \in \R} g_t \a^{ss}(q). \] 

Each leaf $\a^+$ of the strongly unstable foliation $\F^{uu}$, as well
as each leaf $\a^-$ of the strongly stable foliation $\F^{ss}$,
carries a globally defined conditional measure $\mu_{\a^+}$, or
$\mu_{\a^-}$, which is $\G$-invariant, and has the property that
\begin{align*}
(g_t)_* \mu_{\a^+} & = e^{ht} \mu_{g_t \a^+} \\
(g_t)_* \mu_{\a^-} & = e^{-ht} \mu_{g_t \a^-},
\end{align*}
where $h = 6g-6+2b$ is the topological entropy of the the geodesic flow
$g_t$ on $Q / \G$, with respect to $\mu$.  Each leaf of $\F^s$ is
homeomorphic to an open subset of $\MF$ via the map $\eta^-$, so the
pullback of the Thurston measure $\nu$ on $\MF$ defines a conditional
measure on $\F^s$. The foliations $\F^u$ and $F^{ss}$ form a
complementary pair in the sense of Margulis \cite{margulis}, as do
$\F^s$ and $\F^{uu}$.

Let $q \in Q$, and let $\a^u(q)$ be the leaf of the unstable foliation
through $q$. By the Hubbard-Masur Theorem, the projection $\pi:Q \to
\T$ induces a smooth bijection between $\a^u(q)$ and $\T$. The
globally defined conditional measure $\mu_{\a^u(q)}$ on the leaf
$\a^u(q)$ projects onto a measure on the \Teich space, which is
absolutely continuous with respect to the smooth measure $\m$, so we
may consider the Radon-Nikodym derivative of $\pi_*(\mu_{\a^u(q)})$
with respect to $\m$. Let $\l^+ : Q \to \R$ be the function defined by
\[ \frac{1}{\l^+(q)} = \frac{d(\pi_*(\mu_{\a^u(q)}))}{d \m}, \]
and similarly define $\l^-$ to be the function
\[ \frac{1}{\l^-(q)} = \frac{d(\pi_*(\mu_{\a^s(q)}))}{d \m}, \]
where the Radon-Nikodym derivatives are evaluated at $\pi(q)$. Let
$s_x$ be the conditional measure of $\mu$ on $S(x)$. The
\emph{Hubbard-Masur function} $\Lambda$ is defined to be
\[ \Lambda(x) = \int_{S(x)} \l^+(q) \ ds_x(q) = \int_{S(x)} \l^-(q) \ ds_x(q). \]
The functions $\l^+, \l^-$ and $\Lambda$ are all $\G$-invariant.

Athreya, Bufetov, Eskin and Mirzakhani \cite{abem}, showed how to
count the number of images of $x$ in the ball of radius $r$ in \Teich
space, with the \Teich metric.  If $X$ is a finite subset of $\G$, we
will write $\norm{X, \text{\emph{condition}}}$ to denote the number of
elements $\g \in X$ which also satisfy \emph{condition}. Let $B_r(x)$
be the ball of radius $r$ in the \Teich metric, centered at $x$ in
\Teich space, and let $U \subset S(x)$ and $V \subset S(y)$ be Borel
sets. We shall write $\d U$ for the boundary of $U$, which is
$\overline{U} \cap \overline{S(x) \setminus U}$, and we shall always
require that the sets $U$ and $V$ have boundaries of measure zero,
with respect to either $s_x$ or $s_y$, as appropriate.  Given $\g \in
\G$, we will write $q(x, \g y)$ to denote the pair of quadratic
differentials $(q_x(\g y), q_y(\g^{-1} x)) \in S(x) \cross S(y)$, and
so $\g y \in \sect{x}{U}$ if and only if $q(x, y) \in U \cross S(y)$.
The following result is shown in \cite{abem}.

\begin{theorem}[\cite{abem}*{Theorem 7.2}] \label{theorem:abem implicit}
Let $x,y \in \T$, and let $U \subset S(x)$ and $V \subset S(y)$ be
Borel sets whose boundaries have measure zero. Then as $r \to \infty$,

\[ \norm{ \G y \cap B_r(x), \ q(x, \g y) \in U \cross V } \sim
\frac{1}{h} e^{hr} \int_U \l^{-}(q) ds_x(q) \int_V \l^+(q) ds_y(q). \]
\end{theorem}

In \cite{abem} the result is stated for closed surfaces, but the proof
also works for non-sporadic surfaces of finite type.

Veech \cite{veech} showed that the \Teich geodesic flow is
mixing, i.e.  let $\a$ and $\b$ be in $L^2(Q / \G)$. Then
\[ \lim_{t \to \infty} \int_{Q / \G} \a(g_t q) \b(q) d \mu(q) =
\int_{Q / \G} \a(q) d \mu(q) \int_{Q / \G} \b(q) d \mu(q). \]
Following Eskin and McMullen \cite{emc}, we now observe that the
\Teich geodesic flow is also mixing for conditional measures.

\begin{proposition} \label{prop:conditional mixing}
Let $x \in \T$, let $\a$ and $\b$ be continuous non-negative
functions on $Q / \G$, with compact support, and let $U \subset S(x)$
be a Borel set whose boundary has measure zero.  Then
\[ \lim_{t \to \infty} \int_{U} \a(g_t q) \b(q) ds_x (q) =
\int_{Q / \G} \a(q) d \mu(q) \int_{U} \b(q) ds_x(q). \]
\end{proposition}

\begin{proof}
Suppose that $U$ is an open set. Let $U_\e = \bigcup_{\norm{s} \le \e}
g_s U$, and let $I_\e$ be a continuous approximation to the
characteristic function of $U_\e$, i.e. $I_\e$ has maximum value one
at all points of $U_\e$ and is zero outside a small neighbourhood of
$U_\e$. By the definition of conditional measure,
\[ \lim_{t \to \infty} \int_{U} \a(g_t q) \b(q) ds_x (q) =
\lim_{t \to \infty} \lim_{\e \to 0} \frac{1}{2 \e} \int_{Q / \G} \a(g_t
q) I_\e(q) \b(q) d \mu (q), \]
As $\a(q)$ and $I_\e(q)\b(q)$ are continuous functions with compact
support, they are almost constant on sufficiently short segments of
the geodesic flow, i.e. for all $\delta > 0$ there is an $\e > 0$ such
that $ \norm{ \a(g_s q) - \a(q) } \le \delta $ and $ \norm{ I_\e(g_s
q) \b(g_s q) - I_\e(q) \b(q) } \le \delta $ for all $q \in Q/ \G$, and
for all $\norm{s} \le \e$.  As the geodesic flow preserves the lengths
of flow line segments, this implies that $\norm{ \a(g_{t+s} q) -
\a(g_t q) } \le \delta$ for all $q \in Q / \G$, $\norm{s} \le \e$, and
for all $t$.  Therefore the inner limit convergences uniformly,
independently of $t$, and so we may swap the order of the limits. The
mixing property of the \Teich flow implies
\[ \lim_{\e \to 0} \lim_{t \to \infty} \frac{1}{2 \e} \int_{Q / \G}
\a(g_t q) I_\e(q) \b(q) d \mu (q) = \lim_{\e \to 0} \frac{1}{2 \e}
\int_{Q / \G} \a(g_t q) d \mu(q) \int_{Q / \G} I_\e(q) \b(q) d \mu (q), \]
and by the definition of conditional measure, this is equal to
\[ \int_{Q / \G} \a(q) d \mu(q) \int_{U} \b(q) ds_x(q). \]
The result now follows for Borel sets with measure zero boundaries by
approximating them by open sets.
\end{proof}

\section{Limit sets} \label{section:limit sets}

Let $R$ be a set of elements of $\G$ for which there is an upper bound
on their translation distance on the complex of curves.  Let $R_k y$
be the $k$-dense subset of $R y $ in \Teich space, i.e. $R_k$ consists
of those $\g \in R$ such that there is some other element $\g' \in R$
with $\dt{\g y, \g' y} \le k$. We will write $R'_k y$ for the
complement of $R_k y$ in $R y$, and so $R'_k y$ is a $k$-separated
subset of $\T$, as any two elements of $R'_k y$ are distance at least
$k$ apart.  The main aim of this section is to show that the limit set
of $R_k y$ in the visual boundary $S(x)$ has $s_x$-measure zero, Lemma
\ref{lemma:rkzero}, where $s_x$ is the conditional measure induced by
the Masur-Veech measure $\mu$ on $Q$.  We will also show that the
limit set of a metric $k$-neighbourhood of geodesic ray in the visual
boundary is contained in the zero set of the vertical foliation of the
geodesic ray, Lemma \ref{lemma:zero set}. We start by reviewing the
properties of some useful spaces associated to the \mcg.

By work of Masur and Minsky \cite{mm1}, the following three spaces
associated with a surface are quasi-isometric $\delta$-hyperbolic
spaces.

\begin{itemize}

\item The complex of curves $\C(\S)$, which is a simplicial complex
  whose vertices are isotopy classes of simple closed curves, and
  whose simplices are spanned by collections of disjoint simple closed
  curves. We shall write $d_{\C}$ for the metric induced on the
  $1$-skeleton of $\C(\S)$ by setting every edge length equal to one.

\item Electrified \Teich space $\T_{el}$, which is the metric space
  obtained by adding a vertex $v_\a$ for every isotopy class of simple
  closed curve $\a$ on $\S$, and then adding an edge of length
  one half connecting every point of $Thin_\e(\a)$ to $v_\a$.

\item The relative space $\Grel$, which consists of the \mcg $G$, with
  a word metric $\widehat{d}$ coming from the union of a finite
  generating set for $G$, together with a collection of subgroups
  consisting of stabilizers of representatives of simple closed curves
  under the action of the \mcg. We will refer to $\dhat{1, \g}$ as the
  relative length of $\g$.

\end{itemize}

Klarreich \cite{klarreich} identified the Gromov boundary of these
spaces, which we now describe.  Thurston constructed a boundary for
$\T$ consisting of the space of projective measured foliations, which
we shall denote $\PMF$. The space $\PMF$ is a sphere of dimension
$6g-7+2b$, and $\T \cup \PMF$ is homeomorphic to a ball on which $\G$
acts continuously. The Thurston and visual boundaries are distinct,
see for example Kerckhoff \cite{kerckhoff} and Masur \cite{masur2},
and in particular, the action of $\G$ does not extend continuously to
the visual boundary.  There is an inclusion map from \Teich space $\T$
to electrified \Teich space, $\T_{el}$.  This inclusion map does not
extend continuously to the entire Thurston boundary $\PMF$, but
Klarreich \cite{klarreich} shows that it does extend continuously to
the set of minimal foliations $\fmin$ in $\PMF$, and that $\fmin$ is
the Gromov boundary for $\C(\S)$, and hence for the spaces listed
above which are quasi-isometric to $\C(\S)$. The set $\fmin$ is the set of
minimal foliations in $\PMF$, i.e. those foliations for which no
simple closed curve is contained in a (possibly singular) leaf of the
foliation, and furthermore, two foliations are identified if they are
topologically equivalent, even if they have different measures.

Let $Y$ be a subset of the relative space $\Grel$, and let $L$ be a
real number. We define a relative $L$-horoball neighbourhood of $Y$,
which we shall denote $\horo{Y}{L}$, to be the union of balls in
$\Grel$ centered at $y \in Y$, of radius $\nhat{y}+L$, i.e.
\[ \horo{Y}{L} = \bigcup_{y \in Y} \widehat B_{\nhat{y}+L}(y). \] 
We emphasize that this definition uses relative distance in $\Grel$.
We now show that the limit set of $\horo{Y}{L}$ in the Gromov boundary
$\fmin$ is the same as the limit set of $Y$ in $\fmin$.

\begin{lemma}
Let $Y$ be a subset of the relative space $\Grel$. Then the limit set
of $\horo{Y}{L}$ in the Gromov boundary is equal to the limit set of $Y$.
\end{lemma}

\begin{proof}
We will choose the identity element $1$ to be a basepoint in the
relative space $\Grel$, and we will write $\delta$ for the constant of
hyperbolicity of $\Grel$. We will write $\gp{x}{y}$ for the Gromov
product $\tfrac{1}{2}(\dhat{1,x} + \dhat{1,y} - \dhat{x,y})$, which is
equal to the distance from $1$ to a geodesic $[x,y]$, up to an
additive error which only depends on $\delta$, see for example Bridson
and Haefliger \cite{bh}*{Chapter III.H}. Let $x_n$ be a sequence in
$\horo{Y}{L}$ which converges to a point in the Gromov boundary, so in
particular $\dhat{1, x_n}$ tends to infinity. By the definition of
$\horo{Y}{L}$, each $x_n$ lies in a ball of radius $\dhat{x_0, y_n} +
L$, for some $y_n$ in $Y$. As $\dhat{1, x_n}$ tends to infinity, this
implies that $\dhat{1, y_n}$ also tends to infinity. We now show that
the sequence $y_n$ converges to the same limit point as the sequence
$x_n$, by showing that the Gromov product $\gp{x_n}{y_n}$ tends to
infinity.  Suppose not, then possibly after passing to a subsequence,
there is a number $K$ such that $\gp{x_n}{y_n} < K$, for all $n$. Let
$p_n$ be the nearest point projection of $x_n$ to a geodesic
$[1,y_n]$, then it is well known that in a $\delta$-hyperbolic space,
$\dhat{1,p_n}$ is equal to the distance from $1$ to a geodesic $[x_n,
y_n]$, up to an additive error that only depends on $\delta$, see for
example \cite{maher2}*{Proposition 3.2}. Therefore $\dhat{1,p_n}$ is
equal to the Gromov product $\gp{x_n}{y_n}$, and so is at most $K$, up
to an additive error which depends only on $\delta$. This implies
that $\dhat{p_n,y_n}$ is roughly $\dhat{1,y_n} - K$, and as
$\dhat{x_n, y_n} \le \dhat{1, y_n} + L$, and $p_n$ is a nearest point
projection of $x_n$ to $[1, y_n]$, this implies that $\dhat{x_n, p_n}$
is equal to $L+K$, up to additive error that depends only on $\delta$.
Therefore $\dhat{1, x_n}$ is at most $L + 2K$, again up to additive
error that depends only on $\delta$, which contradicts the fact that
$\dhat{1, x_n}$ tends to infinity.
\end{proof}

We now show that the limit set of the orbit of a point in \Teich space
under a horoball neighbourhood of a centralizer has measure zero in
the visual boundary.

\begin{lemma} \label{lemma:measure zero}
Let $\g$ be an element of the \mcg $\G$ which does not lie in the
center of $\G$. Let $H = \horo{C(\g)}{L}$ be a relative horoball neighbourhood of the
centralizer $C(\g)$ in $\Grel$. Then the limit set of $Hy$ in the
visual boundary $\overline \T_x$ has $s_x$-measure zero.
\end{lemma}

\begin{proof}
The set of quadratic differentials with uniquely ergodic initial
foliations has full measure in $S(x)$ with respect to $s_x$, as shown
by Masur \cite{masur} and Veech \cite{veech2}. Therefore, it suffices
to consider the subset of $\overline{H y}$ consisting of limit points
which are uniquely ergodic. It is well known that if a sequence of
points $x_n$ in $\T$ converges to a uniquely ergodic foliation $F \in
\PMF$, then the corresponding sequence of initial measured foliations
also converges to the same uniquely ergodic foliation, see for example
Klarreich \cite{klarreich}*{Proposition 5.3}. Therefore the uniquely
ergodic limit points of $\overline{H y}$ are precisely the uniquely
ergodic foliations in the limit set of $H$ in the Gromov boundary
$\fmin$. This in turn is equal to the limit set of $C(\g)$ in the
Gromov boundary, which is contained in the fixed set of $\g$ by
\cite{maher1}*{Proposition 2.5}. The conditional measure $s_x$ on
$S(x)$ is induced from the measure $\mu$ on $Q$, which in turn is
defined in terms of the Thurston measure $\nu$ on $\MF$, and so in order
to show a subset $U$ of $S(x)$ has $s_x$-measure zero, it suffices to show
that the corresponding cone over the uniquely ergodic foliations in
$U$ has $\nu$-measure zero in $\MF$.  If $\g$ is \pA, then the fixed
set consists of two points, which has measure zero. If $\g$ is
reducible, then the fixed set contains no uniquely ergodic foliations,
and so has measure zero.  Finally, if $\g$ is a non-central periodic element, then the fixed
set of $\g$ acting on $\MF$ is a linear subspace of positive
codimension, so has $\nu$-measure zero, and hence the fixed set of
$\g$ in $S(x)$ has $s_x$-measure zero.
\end{proof}

Elements of the \mcg act as simplicial isometries on the complex of
curves.  Recall that the translation length of an isometry $\g$ is
\[ \t_\g = \lim_{n \to \infty} \tfrac{1}{n} \dc{x, \g^n x}, \] 
and this is independent of the choice of point $x \in C(\S)$.  We now
show that for the \mcg, the translation length of an element $\g$ is
coarsely equivalent to the shortest relative length of any conjugate
of $\g$, which we shall denote $\conj{\g}$.

\begin{proposition} \label{prop:coarse}
There are constants $K_1$ and $K_2$, which depend on $\S$, such that for
any element $\g$ in the mapping class group 
\[ \frac{1}{K_1} \t_g - K_2 \le \conj{\g} \le K_1 \t_\g + K_2, \]
where $\t_\g$ is the translation length of $\g$ on the complex of
curves, and $\conj{\g}$ is the shortest relative length of any conjugate of $\g$.
\end{proposition}

\begin{proof}
The quasi-isometry from $\Grel$ to $\C(\S)$ is defined to be the map
which sends $\g$ to $\g x_0$, for some choice of basepoint $x_0$ in
$\C(\S)$, so there are constants $K$ and $k$ such that
\[ \frac{1}{K} \dhat{\g, \g'} -k \le \dc{\g x_0, \g' x_0} \le K
\dhat{\g, \g'} + k, \]
for all $\g$ and $\g'$ in $\G$.

If $\g$ is periodic, then the translation length of $\g$ on $\C(\S)$
is zero.  There are only finitely many conjugacy classes of periodic
elements in the \mcg, so the proposition holds for periodic elements as
long as $K_2$ is at least the maximum of $\conj{\g}$ over all periodic
elements $\g$.

If $\g$ is reducible, then again its translation length on $\C(\S)$ is
zero. Every reducible element preserves a collection of disjoint
simple closed curves, so there is a simple closed curve which is moved
distance at most one by $\g$. The \mcg $\G$ acts coarsely transitively
on $\C(\S)$, in fact every simple closed curve may be moved to within
distance one of the basepoint $x_0$, so this implies that there is a
conjugate $\g'$ of $\g$ with $\dc{x_0, \g' x_0} \le 3$. This implies
that $\dhat{1, \g'} \le K(3 + k)$, and so $\conj{\g} \le K( 3 + k
)$. Therefore the proposition holds for reducible elements, as long as
$K_2$ is at least $K ( 3 + k )$.

Finally, we consider the case in which $\g$ is \pA. We may assume that
we have chosen $\g$ such that $\dhat{1, \g} = \conj{\g}$.  By the
triangle inequality, the distance $\g$ moves any point in $\C(\S)$ is
an upper bound for the translation length of $\g$, so $\t_\g \le
\dc{x_0, \g x_0} \le K \dhat{1, \g} + k$. This gives the left hand
inequality in the proposition above, with $K_1 = K$, and $K_2 = k$.
By work of Masur and Minsky \cite{mm2}, $\g$ has an axis $\a_\g$,
which is a bi-infinite geodesic such that $\a_\g$ and $\g \a_\g$ are
$2 \delta$-fellow travellers, where $\delta$ is the constant of
hyperbolicity for $\C(\S)$. The \mcg acts coarsely transitively on
$\C(\S)$, so with loss of generality we may assume that we have chosen
$\g$ such that its axis $\a_\g$ passes within distance one of the
basepoint $x_0$. We will Let $p_n$ be a closest point on $\a_{\g}$ to
$\g^n x_0$. As $\a_\g$ is an axis for any power of $\g$, this implies
that $\dc{\g^n x_0, p_n} \le 2 \delta$ for all $n$, and so $\dc{q_n,
  q_{n+1}} \ge \dc{q_0, q_1} - 8 \delta$ for each $n$. As the $q_n$
all lie on a common geodesic, this implies that $\dc{q_0, q_n} \ge
n(\dc{q_0, q_1} - 8 \delta)$, and so $\t_\g \ge \dc{q_0, q_1} - 8
\delta$, and this in turn implies that $\t_\g \ge \dc{x_0, \g x_0} -
10 \delta - 2$.  Therefore $\t_\g \ge \tfrac{1}{K} \dhat{1, \g} - k -
10 \delta - 2$.  This gives the right hand inequality in the
proposition, with $K_1 = K$, and $K_2 = K(k + 10 \delta + 2)$.

We have shown that the inequalities in the proposition hold for each
of the three types of elements of the \mcg, so if we choose the
constants $K_1$ and $K_2$ to be the maximum of the constants we have
obtained in each of the three cases above, then the inequalities hold
for all elements of the \mcg, as required.
\end{proof}

Proposition \ref{prop:coarse} above shows that a set $R$ of \mcg
elements which have an upper bound on their translation distance on
$\C(\S)$, is also a set of elements for which there is an upper bound
on the shortest relative length of elements in the conjugacy class of
each element of $R$. Theorem \ref{theorem:dist} below gives
information about the distribution of such elements inside the \mcg.

\begin{theorem} \cite{maher1}*{Theorem 4.1} \label{theorem:dist}
Let $R$ be a set of elements in the \mcg $\G$, each of which is conjugate to
an element of relative length at most $B$. Then, given $B$ and $\g \in \G$,
there is a constant $L$ such that  $R \cap R \g$ is contained in an
$L$-horoball neighbourhood of the centralizer of $\g$.
\end{theorem} 

We now use this to show that the limit set in the visual boundary
$S(x)$ of the $k$-dense elements $R_k y$ in $\T$ has $s_x$-measure
zero.

\begin{lemma} \label{lemma:rkzero}
For a non-sporadic surface whose \mcg has trivial center, 
the closure of $R_k y$ in $\overline \T_x$ has $s_x$-measure zero in
$S(x)$.
\end{lemma}

\begin{proof}
Let $\rho \in R_k$, then, by the definition of $R_k$, there is a
$\rho' \in R \setminus \rho$ such that $\dt{\rho y, \rho' y} \le k$,
so $\rho = \rho' \g$ for some non-trivial $\g$ with $\dt{y, \g y} \le k$. There
are only finitely many such $\g$, as there are only finitely many
elements of $\G y$ in $B_k(y)$.  Therefore $R_k$ is contained in a
finite union of sets of the form $R \cap R \g$, for elements $\g$ with
$\dt{y, \g y} \le k$. By Theorem \ref{theorem:dist}, a set of the form
$R \cap R \g$ is contained in a relative horoball neighbourhood of the
centralizer of $\g$, i.e.  $\horo{C(\g)}{L}$, where the constant $L$
depends on $R$, $\dhat{1, \g}$ and $\G$. As no element $\g$ is
central, Lemma \ref{lemma:measure zero} implies that the limit set of
$(R \cap R \g) y$ has $s_x$-measure zero in $S(x)$, and as the limit
set of $R_k y$ is contained in a finite union of these limit sets, the
limit set of $R_k y$ also has $s_x$-measure zero.
\end{proof}

Finally, it is well known that the limit set in the visual boundary
$S(x)$ of a $k$-neighbourhood of a geodesic ray based at $x$ is
contained in the zero set of the vertical foliation of the geodesic
ray. We will write $Z(F)$ for the set of foliations with zero
intersection number with $F$, and given $X \subset \T$ we will write
$N_k(X)$ to be a metric $k$-neighbourhood of $X$, using the \Teich
metric.

\begin{lemma} \label{lemma:zero set}
Let $q_t$ be a geodesic ray based at $x$ with vertical foliation
$F = \etabar^+(q_t)$. Then the limit set of $N_k(q_t)$ in the visual
boundary $S(x)$ is contained in $Z(F)$.
\end{lemma}

\begin{proof}
This follows immediately from work of Ivanov \cite{ivanov}, which
shows that if $p$ and $q$ are quadratic differentials with vertical
foliations with non-zero intersection number, then the corresponding
\Teich geodesic rays $p_t$ and $q_t$ diverge, i.e. $\dt{p_t, q_t} \to
\infty$, as $t \to \infty$. In fact, divergent \Teich rays are
completely classified in terms of their vertical foliations, see
Ivanov \cite{ivanov}, Lenzhen and Masur \cite{lm} and Masur
\cites{masur2, masur3}.
\end{proof}

\section{Asymptotics} \label{section:asymptotics}

Let $R$ be a subset of $\G$ consisting of elements for which there is
an upper bound on their translation distance on the complex of curves.
In this section we will prove Theorem \ref{theorem:main}, i.e. we will
show that that the proportion of lattice points $\G y$ in $B_r(x)$
which lie in $R y$ tends to zero as $r$ tends to infinity.

As before, let $R_k$ be the subset of $R$ consisting of those elements
whose images in $\T$ are $k$-dense in $\T$, and let $R'_k y$ be the
complement $R y \setminus R_k y$, and so $R'_k y$ is a $k$-separated
subset of $\T$. In Lemma \ref{lemma:rkzero} we showed that the limit
set of $R_k y$ in the visual boundary $S(x)$ has $s_x$-measure zero.
We now show that this implies that the proportion of lattice points in
$\G y \cap B_r(x)$ which lie in $R_k y$ tends to zero as $r$ tends to
infinity.

\begin{lemma} \label{lemma:closed}
Let $x,y \in \T$, and let $U \subset S(x)$ and $V \subset S(y)$ be
Borel sets whose boundaries have measure zero. Let $X$ be a closed
subset of $\overline \T_x$ such that $X \cap S(x)$ has measure zero
with respect to $s_x$. Then
\[ \frac{ \norm{X \cap \Gamma y \cap B_r(x), \ q(x, \g y) \in U \cross V
  }}{\norm{\Gamma y \cap B_r(x), \ q(x, \g y) \in U \cross V }} \to 0,
\text{ as } r \to \infty. \]
\end{lemma}

\begin{proof}
Consider an open set $W \subset U$ whose closure is disjoint from $X$.
As $X$ is closed, there is a number $D$, depending on $V$, such that
$\sect{x}{W} \setminus B_D(x)$ is disjoint from $X$. Then the
proportion of lattice points inside $B_r(x) \cap \sect{x}{W}$ which
lie in $X$ is at most the proportion of lattice points in $B_r(x) \cap
\sect{x}{U}$ which lie outside $B_D(x)$, and this decays
asymptotically exponentially in $r$, and in particular tends to zero.
Therefore, by Theorem \ref{theorem:abem implicit}, the limiting
proportion of lattice points which lie in $X$ inside $B_r(x) \cap
\sect{x}{U}$ is at most
\[ \lim_{r \to \infty} \frac{ \norm{X \cap \Gamma y \cap B_r(x), \
    q(x, \g y) \in U \cross V }}{\norm{\Gamma y \cap B_r(x), \ q(x, \g
    y) \in U \cross V }} \le \frac{ \int_{U \setminus W}
  \l^-(q)ds_x(q) } { \int_{U} \l^-(q)ds_x(q) } . \]
The function $\l^-$ is absolutely continuous with respect to $s_x$, and we may
choose a sequence of open sets $W_n \subset U$ such that the
$s_x$-measure of $W_n$ tends to the $s_x$-measure of $U$. This implies
that the proportion of lattice points in $X$ tends to zero, as $r$
tends to infinity.
\end{proof}

We now complete the proof of Theorem \ref{theorem:main}.

\begin{proof}
We will first assume that the mapping class group has trivial center,
which covers all cases except for  the genus two
surface and the twice-punctured sphere, which we consider at the very end.
Let $R_k y$ be the $k$-dense subset of $R y$ in \Teich space $\T$, and
let $R'_k y$ be its complement $R y \setminus R_k y$, which is a
$k$-separated set in $\T$. Therefore we may rewrite the fraction in line
\eqref{eq:main} of the statement of Theorem \ref{theorem:main} as
\[ \frac{ \norm{ R_k y \cap B_r(x), \ q(x, \g y) \in U \cross V } }{
  \norm{ \G y \cap B_r(x), \ q(x, \g y) \in U \cross V } } + \frac{
  \norm{ R'_k y \cap B_r(x), \ q(x, \g y) \in U \cross V } }{ \norm{
    \G y \cap B_r(x), \ q(x, \g y) \in U \cross V } }. \]
We have shown in Lemma \ref{lemma:rkzero} that $\overline{R_k y}$ has
$s_x$-measure zero in $S(x)$, so the first term tends to zero as $r$
tends to infinity for any $k$, by Lemma \ref{lemma:closed}.  Therefore
\begin{equation} \label{eq:split} 
\lim_{r \to \infty } \frac{ \norm{ R y
      \cap B_r(x), \ q(x, \g y) \in U \cross V } }{ \norm{ \G y \cap
      B_r(x), \ q(x, \g y) \in U \cross V } } = \lim_{r \to \infty}
  \frac{ \norm{ R'_k y \cap B_r(x), \ q(x, \g y) \in U \cross V } }{
    \norm{ \G y \cap B_r(x), \ q(x, \g y) \in U \cross V } }.
\end{equation}
Consider the subset of \Teich space consisting of $B_r(x) \cap
\sect{x}{U}$. We divide this set into three parts, which we now
describe.

\begin{itemize}

\item The vertical region: $V_r = N_k(\sect{x}{\d U}) \cap B_r(x)$, where
  $N_k(X)$ is a metric $k$-neighbourhood of $X \subset \T$.

\item The interior region: $I_r = ( \sect{x}{U} \setminus V_r ) \cap B_{r-k}(x) $

\item The annular region: $A_r = \sect{x}{U} \setminus (I_r \cup V_r)$

\end{itemize}

\begin{figure}[H] 
\begin{center}
\epsfig{file=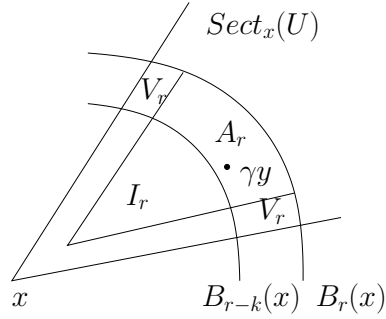, height=120pt}
\end{center}
\caption{Dividing a sector into three regions.}\label{picture1}
\end{figure}

Therefore we may rewrite line \eqref{eq:split} as a sum of three terms.
\begin{equation} \label{eq:three terms} \lim_{r \to \infty} ( \frac{
    \norm{R'_k y \cap V_r, \ q_y(\g^{-1}) \in V} }{ \norm{ \G y \cap
      B_r(x), \ q(x, \g y) \in U \cross V } } + \frac{ \norm{R'_k y
      \cap I_r, \ q_y(\g^{-1}) \in V} }{ \norm{ \G y \cap B_r(x), \
      q(x, \g y) \in U \cross V } } + \frac{ \norm{R'_k y \cap A_r, \
      q_y(\g^{-1}) \in V} }{ \norm{ \G y \cap B_r(x), \ q(x, \g y)
    \in U \cross V } } )
\end{equation}
We shall consider each term in turn. In each case, we find upper
bounds for the term which may depend on $k$, or additional parameters,
and we then show that there is some sequence of upper bounds which
tends to zero. We start by considering the lattice points $R'_k y$ in
the vertical region $V_r$.

\begin{claim} For any fixed $k$, the proportion of lattice
  points in $B_r(x)$ with $q(x, \g y) \in U \cross V$, which lie in both $R'_k y$ and
  the vertical region $V_r$, tends to zero as $r \to \infty$, i.e.
\[ \lim_{r \to \infty} \frac{ \norm{R'_k y \cap V_r, \ q(x, \g y)
    \in U \cross V}}{\norm{ \G y \cap
B_r(x), \ q(x, \g y) \in U \cross V }} = 0. \]
\end{claim}

\begin{proof}
The number of elements of $R'_k y$ in the vertical boundary $V_r$ is at most
the total number of lattice points $\G y$ in $V_r$, i.e.
\[ \norm{ R'_k y \cap V_r, \ q(x, \g y) \in U \cross V } \le
  \norm{ \G y \cap V_r, \ q(x, \g y) \in U \cross V }. \]
The set $V_r$ is contained in $N_k(\sect{x}{\d U})$, and by Lemma
\ref{lemma:zero set}, the limit set of $N_k(\sect{x}{\d U})$ is
contained in $Z(\d U)$, where $Z(\d U)$ is all foliations with zero
intersection with some $F$ in $\d U$, the frontier of $U$. As $\d U$
has $s_x$-measure zero, the set $Z(\d U)$ also has $s_x$-measure
zero. Let $W_1 \supset W_2 \supset \cdots$ be a nested sequence of
open neighbourhoods of $Z(\d U)$, such that $\bigcap W_n = Z(\d
U)$. Then for each $n$ there is a $D_n$ such that
\[ N_k(\sect{x}{\d U}) \setminus B_{D_n}(x) \subset \sect{x}{W_n}, \]
and $B_{D_n}(x)$ contains only finitely many lattice points.  This
implies
\[ \lim_{r \to \infty} \frac{ \norm{ R'_k y \cap V_r, \ q(x, \g y) \in
    U \cross V} }{ \norm{ \G y \cap B_r(x), \ q(x, \g y) \in U \cross
    V } } \le \lim_{r \to \infty} \frac{ \norm{ \G y \cap B_r(x), \
    q(x, \g y ) \in W_n \cross V } }{ \norm{ \G y \cap B_r(x), \ q(x,
    \g y) \in U \cross V }}, \]
for any fixed $n$. Using Theorem \ref{theorem:abem implicit}, we may
take the limit as $r$ tends to infinity on the right hand side, to obtain
\[ \lim_{r \to \infty} \frac{ \norm{ R'_k y \cap V_r, \ q(x, \g y) \in U
    \cross V} }{ \norm{ \G y \cap B_r(x), \ q(x, \g y) \in U \cross V
  } } \le \frac{ \int_{W_n} \l^-(q) ds_x(q) }{ \int_U \l^-(q) ds_x(q)
}, \]
and this upper bound holds for any fixed $n$. However, the top
integral on the right hand side tends to zero as $n$ tends to
infinity, as the $s_x$-measure of $W_n$ tends to zero and $\l^-$ is
absolutely continuous with respect to $s_x$.
\end{proof}

We have shown that the first term on the right hand side of
\eqref{eq:three terms} tends to zero as $r \to \infty$.  We now
consider the lattice points $R'_k y$ in the interior region $I_r$.  We
find an upper bound for the proportion of lattice points $\g y$ in
$B_r(x)$ with $q(x, \g y) \in U \cross V$, which lie in both $R'_k y$
and the interior region $I_r$. This upper bound depends on $k$, and we
show that in fact it decays exponentially in $k$.

\begin{claim} 
There is an upper bound for the limiting proportion of lattice points
in $I_r$, with $q(x, \g y) \in U \cross V$ and which are contained in $R'_k
y$, and furthermore, this upper bound decays exponentially in $k$. In
particular, 
  \[ \lim_{r \to \infty} \frac{ \norm{R'_k y \cap I_r, \ q(x, \g y)
      \in U \cross V} }{\norm{ \G y \cap B_r(x), \ q(x, \g y) \in U
      \cross V }} \to 0, \text{ as } k \to \infty. \]
\end{claim}

\begin{proof}
Every element of $R'_k y$ in the interior set $I_r$ is surrounded by a ball
of radius $k$, which contains $\norm{ \G y \cap B_k(y) }$ lattice
points, at most one of which lies in $R'_k y$. Every ball of radius
$k$ in \Teich space which intersects $I_r$ is contained in $B_r(x)
\cap \sect{x}{U}$. The terminal quadratic differentials for these
lattice points need not be contained in $V$, so we obtain an upper
bound in terms of $U \cross S(y)$ instead of $U \cross V$, namely,
\[ \norm{R'_k y \cap I_r, \ q(x, \g y) \in U \cross V} \le \frac{1}{
  \norm{ \G y \cap B_k(y) } } \norm{\G y \cap B_r(x), \ q(x, \g y) \in
  U \cross S(y)}. \]
Therefore, by Theorem \ref{theorem:abem implicit}, we obtained the
following upper bound, 
\[ \lim_{r \to \infty } \frac{ \norm{R'_k y \cap I_r, \ q(x, \g y) \in
    U \cross V} }{ \norm{\G y \cap B_r(x), \ q(x, \g y) \in U \cross V
  } } \le \frac{ 1 }{ \norm{\G y \cap B_k(y)} } \frac{\int_{S(y)}
  \l^+(q) ds_y(q)}{\int_V \l^+(q) ds_y(q)}. \]
The result now follows by applying Theorem \ref{theorem:abem
  implicit} to $\norm{\G y \cap B_k(y)}$, which shows that the
denominator of the right hand side grows exponentially in $k$.
\end{proof}

Finally, we consider the lattice points in the annular region
$A_r$.

\begin{claim} The limiting proportion of lattice points in the annular region
  $A_r$, with $q(x, \g y) \in U \cross V$ and which lie in $R'_k y$,
  has an upper bound which depends on $k$, and this upper bound tends
  to zero as $k \to \infty$, i.e.
  \[ \lim_{r \to \infty} \frac{ \norm{R'_k y \cap A_r, \ q(x, \g y)
      \in U \cross V} }{\norm{ \G y \cap B_r(x), \ q(x, \g y) \in U
      \cross V }} \to 0, \text{ as } k \to \infty. \]
\end{claim}

\begin{proof}
We will construct a collection of upper bounds $C(k,d)$ for the term
above, which depend on $k$, and an extra parameter $d$, which may be any
positive number less than $k/2$. We will then show that there is a sequence of
these upper bounds which tends to zero, by finding an upper bound for
$\lim_{k \to \infty} C(k,d)$ which depends on $d$, and then showing
that $\lim_{d \to \infty} (\lim_{k \to \infty} C(k,d)) = 0$.

We now give a brief overview of the argument, before giving the
details.  Each lattice point $\g y$ in $R'_k y \cap A_r$ is surrounded
by a ball of radius $k$ disjoint from any other element of $R'_k y$,
and which contains $\norm{\G y \cap B_k(y)}$ lattice points of $\G y$.
However, many of these points may lie outside $B_r(x)$, so we cannot
use $1/\norm{\G y \cap B_k(y)}$ as an estimate for the proportion of
lattice points in $\G y \cap A_r$ which also lie in $R'_k y$.  Recall
that $q_{\g y}(x)$ is the terminal quadratic differential of the
geodesic from $x$ to $\g y$, so $\pi g_{k/2} q_{\g y}(x)$, is the
point distance $k/2$ from $\g y$ along the geodesic from $x$ to $\g
y$. The ball of radius $k/2$ centered at $\pi g_{k/2} q_{\g y}(x)$ is
contained in $B_k(\g y) \cap B_r(x)$, but $\pi g_{k/2} q_{\g y}(x)$ is
not necessarily an element of $\G y$, and a priori, this ball need not
contain any lattice points.  Let $N_d(\G y)$ be a metric
$d$-neighbourhood of the \Teich lattice. If $\pi g_{k/2} q_{\g y}(x)
\in N_d(\G y)$ then there is a lattice point $\g' y$ distance at most
$d$ from $\pi g_{k/2} q_{\g y}(x)$, and so there are at least $\norm{
  \G y \cap B_{k/2-d}(y)}$ lattice points of $\G y$ in $B_k(\g y) \cap
B_r(x)$, as illustrated below in Figure \ref{picture3}.

\begin{figure}[H] 
\begin{center}
\epsfig{file=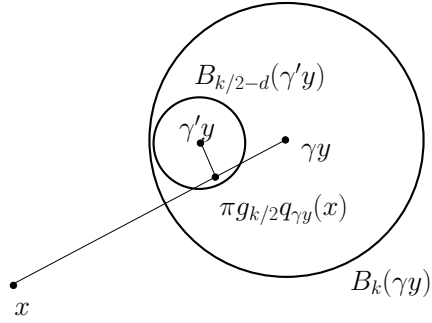, height=120pt}
\end{center}
\caption{The case in which $\pi g_{k/2} q_{\g y}(x)$ lies in $N_d(\G y)$.}\label{picture3}
\end{figure}

Therefore, for lattice points $\g y \in R'_k y$ with $\pi g_{k/2}
q_{\g y}(x) \in N_d(\G y)$, we may use $1/\norm{\G y \cap
  B_{k/2-d}(y)}$ as an upper bound for the proportion of lattice
points $\G y$ in $A_r$ which lie in $R'_k y$.  The number of lattice
points in $R'_k y \cap A_r$ with $\pi g_{k/2} q_{\g y}(x) \not \in
N_d(\G y)$ is bounded above by the number of lattice points in $\G y
\cap A_r$ with $\pi g_{k/2} q_{\g y}(x) \not \in N_d(\G y)$, and we
now describe how to estimate the proportion of lattice points with
this property.  If we project down to moduli space $Q / \G$, then the
geodesic segments from $x$ to $\g y$ all arrive at the same point in
$Q / \G$, which we shall call $y$.  Theorem \ref{theorem:abem
  implicit} implies that the limiting distribution of terminal
quadratic differentials $\g^{-1} q_{\g y}(x) = q_{y}(\g^{-1} x)$ in
$S(y)$ is given by $\l^+(q) ds_y(q)$, up to rescaling. The proportion
of these quadratic differentials for which $\pi g_{k/2} q_y(\g^{-1}
x)$ lies in $N_d(\G)$ is then described as the integral of
$I_d(g_{k/2}q)$ over $S(y)$, where $I_d$ is the characteristic
function for the pre-image of $N_d(\G)$ in $Q$, and so we may apply
Proposition \ref{prop:conditional mixing} to take the limit of this
integral as $k$ tends to infinity. The limit is equal to the volume of
the pre-image of $N_d(\G)$ in moduli space $Q / \G$, and this tends to
one as $d$ tends to infinity. This means that the proportion of
$k$-separated lattice points $R'_k y$ for which there are at least
$\norm{\G y \cap B_{k/2-d}(y)}$ nearby lattice points in $\G y \cap
B_r(x)$, but not in $R'_k y$, tends to one as $k$ tends to infinity,
for appropriate choices of $d$ for each $k$.  Therefore there is a
sequence of upper bounds which tend to zero. We now give a detailed
version of this argument.

Recall that $q_x(\g y)$ is the unit area initial quadratic
differential on $x$ for the \Teich geodesic from $x$ to $\g y$, and
$q_{\g y}(x)$ is the corresponding unit area terminal quadratic
differential on $\g y$.  Recall that $\pi$ is the projection map $\pi
: Q \to \T$, and we denote the \Teich geodesic flow on $Q$ by $g_t$,
so $\pi g_t q_{\g y}(x)$ is the point distance $t$ along the geodesic
ray starting at $\g y$ which passes through $x$.  Let $d < k/2$, and
let $N_d(\G y)$ be a metric $d$-neighbourhood of the \Teich lattice
$\G y$. If $\g y \in R'_k y \cap A_r$, and $\pi g_{k/2} q_{\g y}(x)
\in N_d(\G y)$, then there is a lattice point $\g' y$ with $\dt{\pi
  g_{k/2} q_{\g y}(x), \g' y } \le d$. In particular,
\[ B_{k/2 -d}(\g' y) \subset B_k(\g y) \cap B_r(x) \cap
\sect{x}{U}, \] 
as illustrated above in Figure \ref{picture3}.  There are $\norm{\G y
  \cap B_{k/2 - d}(y)}$ lattice points of $\G y$ in $B_{k/2 - d}(\g'
y)$, and at most one of these lattice points lies in $R'_k y$.  

Let $Y(k, d) = \{ q \in S(y) \mid \pi g_{k/2} q \in N_d(\G y) \}$. We
may divide the points of $R'_k y \cap A_r$ into two sets, depending on
whether or not the corresponding terminal quadratic differential
$q_{y}(\g^{-1} x)$ lies in $Y(k, d)$. Therefore $\norm{ R'_k y \cap
  A_r, \ q(x, \g y) \in U \cross V }$ is equal to 
\begin{equation} \label{eq:Y} 
\norm{ R'_k y \cap A_r, \ q(x, \g y) \in U \cross V \cap Y(k, d)
} + \norm{ R'_k y \cap A_r, \ q(x, \g y) \in U \cross V \setminus
  Y(k, d) }.
\end{equation}
We first consider the first term from line \eqref{eq:Y}. For each $\g
y \in R'_k y \cap A_r$, with $q_{y}(\g^{-1} x) \in Y(k, d)$, we get at
least $\norm{ \G y \cap B_{k/2 - d}(y)}$ lattice points of $\G y$
inside $B_k(\g y) \cap B_r(x) \cap \sect{x}{U}$, at most one of which
lies in $R'_k y$. The terminal quadratic differentials for these
lattice points need not lie in $V$, so we obtain an upper bound in
terms of $U \cross S(y)$ instead of $U \cross V$, i.e.
\begin{equation} \label{eq:Y1}
\norm{ R'_k y \cap A_r, \ q(x, \g y) \in U \cross V \cap Y(k,d)}
\le \frac{1}{ \norm{\G y \cap B_{k/2 - d}(y) } } \norm{ \G y \cap
  B_r(x), \ q(x, \g y) \in U \cross S(y) }. 
\end{equation}
For the second term in line \eqref{eq:Y}, we use the fact that $R'_k y
\subset \G y$, and $A_r \subset B_r(x)$, which gives the following
upper bound,
\begin{equation} \label{eq:Y2}
\norm{ R'_k y \cap A_r, \ q(x, \g y) \in U \cross V \setminus
  Y(k, d)} \le \norm{ \G y \cap B_r(x), \ q(x, \g y) \in U \cross
  V \setminus Y(k, d) }. 
\end{equation}
Now applying Theorem \ref{theorem:abem implicit} to lines
\eqref{eq:Y1} and \eqref{eq:Y2} above, and adding them together as in
line \eqref{eq:Y}, we obtain the following upper bound,
\begin{equation} \label{eq:third term}
\lim_{r \to \infty} \frac{ \norm{ R'_k y \cap A_r, \ q(x, \g y) \in
    U \cross V } }{ \norm{ \G y \cap B_r(x), \ q(x, \g y) \in U \cross
    V } } \le \frac{1}{ \norm{\G y \cap B_{k/2 - d}(y) } }
\frac{\int_{S(y)} \l^+(q) ds_y(q)}{\int_V \l^+(q) ds_y(q)} + \frac{
  \int_{ V \setminus Y(k, d) } \l^+(q) ds_y(q) }{ \int_V \l^+(q)
  ds_y(q) },  
\end{equation}
and this upper bound holds for all $k$, for all $d < k/2$.

By Theorem \ref{theorem:abem implicit}, the first term on the right
hand side above tends to zero as $k$ tends to infinity, for any fixed
$d$. We now consider the second term above, which we may rewrite as $1
- p(k, d)$, where 
\begin{equation} \label{eq:pkd}
p(k, d) = \frac{ \int_{ V \cap Y(k, d) } \l^+(q) ds_y(q) }{
  \int_V \l^+(q) ds_y(q) }. 
\end{equation}
Recall that $Y(k, d) = \{ q \in S(y) \mid \pi g_{k/2} q \in N_d(\G y)
\}$, and let $I_d$ be the characteristic function for the pre-image of
$N_d(\G y)$ in $Q$. The function $I_d$ is $\G$-equivariant, so we may
write the expression for $p(k,d)$ in line \eqref{eq:pkd} above as
\[ p(k,d) = \frac{ \int_V I_d(g_{k/2} q) \l^+(q) ds_y(q) }{\int_V
  \l^+(q) ds_y(q)}. \]
Therefore as the \Teich geodesic flow is mixing for conditional
measures, Proposition \ref{prop:conditional mixing}, this implies
\[ \lim_{k \to \infty} p(k,d) = \int_{ Q / \G } I_d(q) d \mu(q). \]
Therefore, if we fix $d$ and let $k$ tend to infinity in line
\eqref{eq:third term} above, we obtain the following upper bound,
\begin{equation} \label{eq:p}
 \lim_{k \to \infty} \lim_{r \to \infty} \frac{ \norm{ R'_k y \cap
    A_r, \ q(x, \g y) \in U \cross V } }{ \norm{ \G y \cap B_r(x), \
    q(x, \g y) \in U \cross V } } \le 1 - \int_{ Q / \G } I_d(q) d
\mu(q), 
\end{equation}
and this upper bound holds for all $d$. However, the integral on line
\eqref{eq:p} above tends to one as $d$ to infinity, as we have
normalized our measures so that the volume of moduli space is one.
Therefore, the term on the right hand side of \eqref{eq:p} above tends
to zero as $d$ tends to infinity. As the right hand side of
\eqref{eq:p} is an upper bound for the left hand side, for all values
of $d$, this implies that the left hand side of \eqref{eq:p} is zero,
as required.
\end{proof}

Therefore, all three terms on the right hand side of \eqref{eq:three
  terms} tend to zero as $k$ tends to infinity. This completes the proof of
Theorem \ref{theorem:main} in the case that the \mcg has trivial
center.

We now deal with the two cases in which the \mcg has non-trivial
center, which are the genus two surface $\S_{2,0}$, and the
twice-punctured torus $\S_{1,2}$, and we will write $\G_{g,b}$ for the
mapping class group of $\S_{g,b}$. In the case of the genus two
surface, the center $Z$ of $\G_{2,0}$ is $\Z / 2 \Z$, generated by the
hyperelliptic involution, and the quotient of the genus two surface by
the hyperelliptic involution is the six-punctured sphere. The
hyperelliptic involution acts trivially on $\T(\S_{2,0})$, which is
isometric to $\T(\S_{0,6})$ and $\G_{2,0} / Z$ is isomorphic to the
\mcg of the six-punctured sphere, so a \Teich lattice in
$\T(\S_{2,0})$ is isometric to a \Teich lattice in $\T(\S_{0,6})$. The
hyperelliptic involution also acts trivially on the complex of curves,
so $\C(\S_{2,0})$ is isometric to $\C(\S_{0,6})$. Therefore a set $R$
of elements of bounded translation length length on $\C(\S_{2,0})$ is
also a set of elements of bounded translation length on
$\C(\S_{0,6})$. Theorem \ref{theorem:main} holds for $\G_{0,6}$, as
$\G_{0,6}$ has trivial center, and so this implies that Theorem
\ref{theorem:main} also holds for $\G_{2,0}$.  In the case of the
twice-punctured torus, the quotient surface under the hyperelliptic
involution is the five-punctured sphere, and again the hyperelliptic
involution acts trivially on \Teich space and the complex of
curves. Therefore the argument above works exactly as before, except
for the fact that $\G_{1,2} / Z$ is a finite index subgroup of
$\G_{0,5}$. However, by Theorem \ref{theorem:abem implicit}, the
asymptotic number of lattice points of $\G_{1,2}y \cap B_r(x)$ in a
bisector is a constant multiple of the asymptotic number of lattice
points of $\G_{0,5}y \cap B_r(x)$ in the same bisector, and so the
proportion of points in $Ry \cap B_r(x)$ in the bisector tends to zero
in either lattice. This completes the proof of Theorem
\ref{theorem:main}.
\end{proof}



\end{document}